

\input amstex

\documentstyle{amsppt}

\magnification=\magstep1

\loadbold

\pageheight{9.0truein}
\pagewidth{6.5truein}


\def\Prim{\operatorname{Prim}}
\def\Spec{\operatorname{Spec}}
\def\ann{\operatorname{ann}}
\def\Max{\operatorname{Max}}
\def\Mod{\operatorname{Mod}}
\def\End{\operatorname{End}}
\def\nxn{n{\times}n}

\def\RqG{R_q[G]}
\def\OqG{\Cal{O}_q(G)}

\def\CqpG{\Bbb{C}_{q,p}[G]}
\def\C{\Cal{C}}
\def\p{\frak{p}}\def\q{\frak{q}}
\def\q{\frak{q}}
\def\a{\frak{a}}
\def\b{\frak{b}}
\def\Art{1}
\def\BroGooone{2}
\def\BroGootwo{3}
\def\Dix{4}
\def\Goo{5}
\def\GooLetone{6}
\def\GooLettwo{7}
\def\GooSta{8}
\def\GooWar{9}
\def\HodLevone{10}
\def\HodLevtwo{11}
\def\HodLevTor{12}
\def\Josone{13}
\def\Jostwo{14}
\def\Mat{15}
\def\McCRob{16}
\def\Pro{17}
\def\Ren{18}
\def\Ros{19}
\def\Rowone{20}
\def\Rowtwo{21}
\def\VdB{22}

\topmatter

\title Noncommutative Images of Commutative Spectra
\endtitle

\abstract We initiate a unified, axiomatic study of noncommutative algebras
$R$ whose prime spectra are, in a natural way, finite unions of commutative
noetherian spectra. Our results illustrate how these commutative spectra can
be functorially ``sewn together'' to form $\Spec R$. In particular, we
construct a bimodule-determined functor $\Mod Z \rightarrow \Mod R$, for a
suitable commutative noetherian ring $Z$, from which there follows a
finite-to-one, continous surjection $\Spec Z \rightarrow \Spec R$. Algebras
satisfying the given axiomatic framework include PI algebras finitely
generated over fields, noetherian PI algebras, enveloping algebras of complex
finite dimensional solvable Lie algebras, standard generic quantum semisimple
Lie groups, quantum affine spaces, quantized Weyl algebras, and standard
generic quantizations of the coordinate ring of $\nxn$ matrices. In all of
these examples (except for the non-finitely-generated noetherian PI algebras),
$Z$ is finitely generated over a field, and the constructed map of spectra
restricts to a surjection $\Max Z \rightarrow \Prim R$.
\endabstract

\rightheadtext{Noncommutative Images}

\author Edward S. Letzter \endauthor

\address Department of Mathematics, Temple University, Philadelphia,
PA 19122 \endaddress

\email letzter\@math.temple.edu \endemail

\thanks The author thanks the Department of Mathematics at the University of
Pennsylvania for its hospitality; this paper was begun while he was a visitor
on sabbatical there in Fall 2004. The author is grateful for support during
this period from a Temple University Research and Study Leave Grant. This
research was also supported in part by a grant from the National Security
Agency. \endthanks

\endtopmatter

\document

\head 1. Introduction \endhead

\subhead 1.1 \endsubhead It is well known that the prime ideal theories of the
following noncommutative algebras have many of the basic properties found in
commutative settings:

\smallskip\itemitem{(a)} noetherian PI algebras and PI algebras
finitely generated over a field (see \S 3 below), 

\smallskip\itemitem{(b)} enveloping algebras of
finite dimensional completely solvable Lie algebras (\S 4), 

\smallskip\itemitem{(c)} standard
generic quantized coordinate rings of semisimple Lie groups (\S 5),  

\smallskip\itemitem{(d)} certain other generic quantized coordinate rings,
including quantized Weyl algebras and the standard quantized coordinate rings
of $\nxn$ matrices (\S 6).
\par\smallskip

\noindent Moreover, a longstanding (essentially achieved) goal for all of
these algebras has been the parameterization of their prime and primitive
ideals using algebraic geometric data.

Our aim in this paper is to present, and begin the study of, an abstract
axiomatic setting that includes all of these examples. In our main result,
given in (2.18), we show that the prime spectra occurring within this
framework are continuous, functorial, finite-to-one images of spectra of
commutative noetherian rings.

\subhead 1.2 \endsubhead Our analysis throughout is based on the deep
structure theories available for the classes of algebras listed in (a), (b),
(c), and (d). These theories encompass decades of work, by legions of
mathematicians. The reader unfamiliar with any of the relevant material is
directed, for example, to the texts \cite{\BroGooone; \Dix; \GooWar; \Josone;
\McCRob; \Rowone}. 

\subhead 1.3 \endsubhead The specific algebras $R$ listed in (1.1) all have
the following in common: Their prime spectra are unions of finitely many
locally closed subsets (in the Jacobson/Zariski topology), each homeomorphic
to the spectrum of some commutative noetherian ring, and each naturally
homeomorphic to a localization of a factor of $R$. Hypotheses A, B, and C (and
stronger variants), presented in \S 2, specify these properties in detail. A
further hypothesis D, under which $R$ satisfies the Nullstellensatz, allows
applications to primitive ideals.

\subhead 1.4 \endsubhead To more precisely but still briefly summarize our
main results, stated in (2.18), suppose that $R$ satisfies the hypotheses A,
B, and C (referred to in the preceding paragraph and exactly specified in \S
2). Then there is a homomorphism from $R$ into an algebra $T$, whose ideals
are completely determined by their intersection with its noetherian center
$Z$, such that
$$\theta : \Spec Z \quad @> \quad \q \; \longmapsto \; \ann\big( (Z/\q)
\otimes _ZT \big)_R \quad >> \quad \Spec R$$
is a continuous, finite-to-one, surjection. In cases (c) and (d) of (1.1),
$\theta$ is bijective. In all of the cases (a), (b), (c), and (d), except for
non-finitely-generated noetherian PI algebras, we can assume both that
$Z$ is finitely generated over the ground field and that $\theta$ restricts to
a surjection from $\Max Z$ onto $\Prim R$. 

Note that $\theta$ is a consequence of the functor
$$\Mod Z \quad @> \quad {\; - \;} \otimes_ZT \quad >> \quad \Mod R.$$
In the approach to noncommutative algebraic geometry developed in \cite{\Ros,
\VdB} (for example), $\Mod R$ is viewed as a category of quasi-coherent
sheaves on a not-explicitly-defined noncommutative scheme, and the preceding
functor is an {\sl affine map\/} from $\Mod Z$ into $\Mod R$. Recall, using a
theorem of Serre, that $\Mod Z$ can be identified with the category of
quasi-coherent sheaves on the commutative scheme $\Spec Z$.

Applications of the above results to the specific classes of algebras listed
in (1.1) are given in sections 3, 4, 5, and 6.

\subhead 1.5 \endsubhead Let $V$ be a topological space that is a (not
necessarily disjoint) union of subspaces $V_1,\ldots,V_t$, each homeomorphic
to the prime spectrum of a commtuative noetherian ring. The coproduct $V_1
\sqcup \cdots \sqcup V_t$ (i.e., the disjoint union equipped with the topology
induced from the individual spaces $V_\ell$) is homeomorphic to the prime
spectrum of a commutative noetherian ring (i.e., the direct product of the
rings corresponding to the $V_\ell$). Also, there is a continuous surjection
from $V_1 \sqcup \cdots \sqcup V_t$ onto $V$. So, the fact that the spectra we
are considering in this paper are continuous images of commutative noetherian
spectra is, in essence, only an elementary afterthought. What is new in our
results are the explicitly described ring homomorphisms $R \rightarrow T$ and
functors $\Mod Z \rightarrow \Mod R$ that produce these continuous
surjections.

\subhead 1.6 \endsubhead The continuous maps presented in (2.18) are not, in
general, topological quotient maps; see (4.3).

\subhead 1.7 \endsubhead If $R$ is the quantized coordinate ring of affine
$n$-space, over an algebraically closed field $\Bbbk$, and $-1$ is not
included in the multiplicative group generated by the quantizing parameters,
then $\Spec R$ is a topological quotient of $\Spec \Bbbk[x_1,\ldots,x_n]$; see
\cite{\GooLetone}. For more general quantizations $\Cal{O}_q(V)$ of coordinate
rings $\Cal{O}(V)$, Goodearl has conjectured that $\Spec \Cal{O}_q(V)$ is a
topological quotient of $\Spec\Cal{O}(V)$; see \cite{\Goo}. If $R$ is the
enveloping algebra of a finite dimensional complex solvable Lie algebra $\frak
g$, then there is a topological quotient map from $\frak g^*$ onto $\Prim R$,
given by the Dixmier map \cite{\Dix; \Mat}.

\subhead 1.8 \endsubhead In (2.16) we note that our hypotheses A, B, C, and D
imply the Dixmier-Moeglin equivalence, following arguments similar to those
found in \cite{\GooLettwo}.

\subhead 1.9 \endsubhead In \cite{\BroGootwo} an axiomatic framework is
developed for certain quantum algebras. Our approach here can be viewed as a
further generalization of this framework.

\subhead Acknowledgement \endsubhead I am happy to acknowledge many
conversations with Ken Goodearl, over several years, on the topics considered
in this paper.

\head 2. The General Theory  \endhead

Throughout this section, $R$ will denote an associative ring with
identity. 

\subhead 2.1 Preliminaries \endsubhead Let $\Lambda$ be a ring, and let $I$ be
an ideal of $\Lambda$.

(i) The center of $\Lambda$ will be denoted $Z(\Lambda)$. The set of elements
of $\Lambda$ regular modulo $I$ will be denoted $\C_\Lambda(I)$.

(ii) Recall that $n \in \Lambda$ is {\sl normal modulo $I$\/} exactly when
$\Lambda.n = n.\Lambda + I$, and that $c \in \Lambda$ is {\sl central modulo
$I$\/} exactly when $a.c = c.a + I$ for all $a \in \Lambda$. When $I = 0$,
elements normal modulo $I$ are {\sl normal\/} (and elements central modulo $I$
are, of course, central). Also recall that if $n$ is a normal element of
$\Lambda$, and if $Q$ is a prime ideal of $\Lambda$, then $n$ is regular
modulo $Q$ if and only if $n \notin Q$. In particular, a nonzero normal
element of a prime ring is regular.

(iii) (Cf\., e.g., \cite{\GooWar, Chapter 10}.) A multiplicatively closed
subset $D$ of $\Lambda$ is a {\sl right denominator set\/} \cite{\GooWar,
p\. 168--169} provided $D$ is a right Ore set \cite{\GooWar, p\. 82} and the
canonical right-Ore-localization map $\Lambda \rightarrow \Lambda D^{-1}$ is
injective \cite{\GooWar, p\. 167}. If $n$ is a normal and regular element of
$\Lambda$, then it is not hard to show that $\{1, n, n^2, \ldots \}$ is a
right (and left) denominator set of $\Lambda$. Localizations at right
denominator sets of regular elements will be referred to as {\sl classical
right quotient rings\/}.

(iv) Let $\varphi\colon \Lambda \rightarrow \Lambda'$ be a ring
homomorphism. For $a \in \Lambda$ and $b \in \Lambda'$, we will use $a.b$ to
denote $\varphi(a).b$, and we will use $b.a$ to denote
$b.\varphi(a)$. Adapting standard terminology to our uses below, we will say
that the homomorphism is {\sl normalizing\/} if $\Lambda'$ is generated as a
left $\Lambda$-module by elements $b$ such that $\Lambda.b = b.\Lambda$, and
we will say that the homomorphism is {\sl centralizing\/} if $\Lambda'$ is
generated as a left $\Lambda$-module by elements $c$ such that $a.c = c.a$ for
all $a \in \Lambda$.

(v) We will use $\Spec \Lambda$, $\Prim \Lambda$, and $\Max \Lambda$ to
denote, respectively, the sets of prime, right primitive, and maximal ideals
of $\Lambda$; throughout this note, each of these sets will be equipped with
the Zariski/Jacobson topology. In particular, the closed subsets of $\Spec
\Lambda$ will have the form
$$V(E) = V_\Lambda(E) = \{ \, P \in \Spec \Lambda \; \mid \; P \supseteq E \,
\},$$
for subsets $E$ of $\Lambda$. Analogous notation will be used for closed
subsets of $\Prim \Lambda$ and $\Max \Lambda$. We will use {\sl primitive\/}
to mean ``right primitive.''

(vi) For an ideal $I$ of $\Lambda$, we will let $\sqrt{I}$ denote the prime
radical of $I$ (i.e., the intersection of all of the prime ideals containing
$I$).

(vii) We will let $\Bbbk$ denote an arbitrary field.

\subhead 2.2 Setup and Notation ($I_\ell$, $J_\ell$, $X_\ell$, $S_\ell$, $S$,
$Y_\ell$, $Y$, $\mu$) \endsubhead We now begin to describe our general
setup. The notation developed here will remain in effect for the rest of this
section.

(i) Choose ideals $I_1,J_1,\ldots,I_t,J_t$ of $R$ such that
$$\Spec R \; = \; X_1 \cup \cdots \cup X_t,$$
where
$$X_\ell = V(I_\ell) \setminus V(J_\ell),$$
for $\ell = 1,\ldots,t$. Note that we are not requiring $X_1,\ldots,X_t$ to be
(pairwise) disjoint. Equip each $X_\ell$ with the subspace topology.

(ii) For $1 \leq \ell \leq t$, let
$$S_\ell = R/I_\ell,$$
and set
$$S = S_1 \times \cdots \times S_t.$$
Let $e_1,\ldots,e_t$ denote the canonical central idempotents
corresponding to this ring product. We obtain a ring
homomorphism
$$\mu \colon R \; @> \quad r \; \longmapsto \; \left(r + I_1, \ldots ,
r + I_t \right) \quad >> S.$$
The standard theory of finite centralizing extensions (cf., e.g.,
\cite{\McCRob, Chapter 10}) ensures that the function
$$\mu^* \colon \Spec S \; @> \quad Q \; \longmapsto \; \mu^{-1}(Q) \quad >>
\Spec R $$
is continuous, with image $\Spec R$, since $S = R.e_1 + \cdots + R.e_t$, and
since $\ker \mu$ is contained within the prime radical of $R$. Continuity and
image can also be directly derived from the given setup.

(iii) Identify $\Spec S$, homeomorphically, with the disjoint union
$$\Spec S_1 \; \sqcup \; \cdots \; \sqcup \; \Spec S_t,$$
equipped with the topology induced from the spaces $\Spec S_\ell$. Note that
each $X_\ell$ is naturally homeomorphic to the complement in $\Spec S_\ell$ of
$$V_{S_\ell}(J_\ell + I_\ell/I_\ell).$$
We will use $Y_\ell$ to denote this complement. Setting $Y = Y_1 \sqcup \cdots
\sqcup Y_t$, and equipping $Y$ with the topology induced from the spaces
$Y_1,\ldots,Y_t$, we may view $Y$ as a subspace of $\Spec S$. (Note that $Y$
need not be homeomorphic to $\Spec R$, even when $X_1,\ldots,X_t$ are
disjoint, because inclusions among prime ideals in $\Spec R$ may not lift to
inclusions among prime ideals in $Y$.)

(iv) Let $P$ be a prime ideal of $R$. Then the preimage of $P$ under $\mu^*$
is exactly the finite (and nonempty) set of kernels in $S$ of the maps,
$$S \; @> \; \text{projection} \; >> \; S_\ell \; = \; R/I_\ell \; @> \;
\text{natural surjection} \; >> \; R/P,$$
for all $\ell$ such that $I_\ell \subseteq P$. Furthermore, it is not hard to
see that the preimage of $P$ under $\mu^*$ intersects $Y$ nontrivially.

(v) Now consider the restriction
$$\mu^* \colon Y \; \longrightarrow \; \Spec R.$$
This function is continuous, by (ii), and surjective, by (iv).

(vi) Note that the restriction of $\mu^*$ in (v) is bijective if
$X_1,\ldots,X_t$ are disjoint.

\subhead 2.3 \endsubhead For a prime ideal $Q$ of $S$, it is easy to see that
$\mu^*(Q)$ is primitive if and only if $Q$ is primitive, and that $\mu^*(Q)$
is maximal if and only if $Q$ is maximal. (These conclusions also follow from
the standard theory of finite centralizing extensions.)  Therefore, all of the
statements in (2.2) remain valid if the word ``prime'' is replaced with the
word ``primitive'' and if ``Spec'' is replaced with ``Prim.'' We can make a
similar substitution with maximal ideals and maximal spectra.

\subhead 2.4 Additional Setup and Notation ($D_\ell$, $T_\ell$, $T$, $\nu$,
$\eta$, $Z_\ell$, $Z$) \endsubhead The following additional notation -- to be
used in conjunction with later hypotheses -- will be retained for the
remainder of this section: Let $D_\ell$ be a right denominator set of
$S_\ell$, for each $\ell$, and let $T_\ell$ be the right Ore localization of
$S_\ell$ at $D_\ell$. Set $T = T_1\times \cdots \times T_t$. (Note that $T$ is
itself an Ore localization of $S$.) Let $\nu\colon S \rightarrow T$ be the
natural product homomorphism, and let $\eta\colon R \rightarrow T$ be the
composition of $\mu$ with $\nu$.

For each $\ell$, let $Z_\ell$ denote the
center of $T_\ell$, and identify $Z_1 \times \cdots \times Z_t$ with the
center $Z$ of $T$.

\subhead 2.5 \endsubhead We now begin to introduce the hypotheses that will
guide our analysis. We will always explicitly indicate when they are in
effect. 

\proclaim{Hypothesis A} $T_\ell$ is right noetherian for all
$\ell$. Consequently, $T$ is right noetherian. \endproclaim

\proclaim{Hypothesis B} $Y_\ell = \{ \, Q \in \Spec S_\ell
\; \mid \; D_\ell \subseteq \C_{S_\ell}(Q) \, \}$, for all $\ell$.
\endproclaim

\subhead 2.6 \endsubhead Assume hypotheses A and B. (i) It follows (e.g.)
from \cite{\GooWar, 10.17, 10.18} that there is a lattice isomorphism
$$\{ \, \text{semiprime ideals of $T_\ell$} \, \} \longrightarrow \{ \,
\text{semiprime ideals $N$ of $S_\ell$ such that $D_\ell \subseteq
\C_{S_\ell}(N)$} \, \},$$
mapping each semiprime ideal of $T_\ell$ to its natural preimage in $S_\ell$.
When we restrict the preceding function to prime ideals, we obtain a
homeomorphism
$$\Spec T_\ell \longrightarrow Y_\ell.$$

Now let $Q$ be a prime ideal of $S_\ell$. Under the preceding map, the
preimage of $Q$ is $Q.T_\ell$ (by, e.g., \cite{\GooWar, 10.18}).  Next, let
$\overline{D_\ell}$ denote the image modulo $Q$ of $D_\ell$. It follows from
hypothesis B that $\overline{D_\ell}$ is a right Ore set, in $S_\ell/Q$, of
regular elements. Hence,
$$T_\ell/Q.T_\ell \; \cong \; (S_\ell/Q) \otimes _{S_\ell}T_\ell \; \cong \;
(S_\ell/Q)\overline{D_\ell}^{-1}.$$
In particular, $T_\ell/Q.T_\ell$ is a classical right quotient ring for
$S_\ell/Q$.

(ii) Now consider the homomorphism $\nu\colon S \rightarrow T$ described in
(2.4). We see that
$$\nu^* \colon \Spec T \; @> \quad \p \; \longmapsto \; \nu^{-1}(\p) \quad >>
\; \Spec S$$
is homeomorphic onto its image, $Y$. Let $Q$ be a prime ideal of $S$. It is
not hard to prove, given (i), that the preimage of $Q$ under $\nu^*$ is $Q.T$
and that $T/Q.T$ is a classical right quotient ring of $S/Q$.

(iii) The necessity in (i) and (ii) of right (or left) noetherianity, as
assumed in hypothesis A, is discussed in \cite{\GooWar, p\. 179}.

\proclaim{2.7 Lemma} Assume hypotheses {\rm A} and {\rm B}. Consider the
function
$$\eta^* \colon \Spec T \; @> \quad \p \; \longmapsto \; \eta^{-1}(\p) \quad
>> \; \Spec R .$$

{\rm (i)} $\eta^* = \mu^* \circ \nu^*$.

{\rm (ii)} $\eta^*$ is a continuous, finite-to-one, surjective function.

{\rm (iii)} When $X_1,\ldots,X_t$ are disjoint, $\eta^*$ is
injective. \endproclaim

\demo{Proof} First, part (i) follows from (2.2) and (2.6ii).

Next, in view of (i), the continuity and surjectivity of $\eta^*$ follow from
(2.6ii) and (2.2v). That $\eta^*$ is finite-to-one follows from (i) and
(2.2iv). Part (ii) follows.

Part (iii) follows from (i) and (2.2vi). \qed\enddemo

\subhead 2.8 \endsubhead As noted in (2.1ii), a normal element $n$ of a ring
$\Lambda$ is regular modulo a prime ideal $Q$ of $\Lambda$ if and only if $n
\notin Q$.  Also, as noted in (2.1iii), the non-negative powers of a regular
normal element form a right denominator set. It follows from these two facts
that the following implies hypothesis B:

\proclaim{Hypothesis B$'$} For each $\ell$, the ideal $I_\ell + J_\ell$ is
generated by a single element $b_\ell$ normal and regular modulo $I_\ell$, and
$D_\ell = \{1+I_\ell,b_\ell+I_\ell,b_\ell^2+I_\ell,\ldots\} \subset S_\ell$.
\endproclaim

\remark{{\bf 2.9} Remarks} (i) When B$'$ holds and $R$ is a finitely generated
$\Bbbk$-algebra, then $T$ is also a finitely generated $\Bbbk$-algebra.

(ii) When B$'$ holds, the homomorphism $\eta\colon R \rightarrow T$ is
normalizing. Moreover, if each of the $b_\ell$ is central modulo $I_\ell$ then
the homomorphism is centralizing.
\endremark

\subhead 2.10 \endsubhead Let $\Lambda$ be a ring, and suppose there is a
product preserving, one-to-one correspondence between the sets of ideals of
$\Lambda$ and $Z(\Lambda)$, given by mutually inverse assignments
$$\a \; \longmapsto \; \a \cap Z(\Lambda), \quad \text{and} \quad \b \;
\longmapsto \; \b \Lambda,$$
for arbitrary ideals $\a$ of $\Lambda$ and $\b$ of $Z(\Lambda)$. We will say
in this situation that $\Lambda$ is {\sl ideally commutative\/}. There is then
a homeomorphism between $\Spec \Lambda$ and $\Spec Z(\Lambda)$, given by
$$\a \; \longmapsto \; \a \cap Z(\Lambda), \quad \text{and} \quad \b \;
\longmapsto \b \Lambda,$$
for prime ideals $\a$ of $\Lambda$ and $\b$ of $Z(\Lambda)$. This
homeomorphism restricts to a homeomorphism between $\Max \Lambda$ and $\Max
Z(\Lambda)$. If, in addition, $Z(\Lambda)$ is a finitely generated
$\Bbbk$-algebra, we will say that $\Lambda$ is {\sl ideally $\Bbbk$-affine
commutative\/}.

We now present the third hypothesis, and a stronger variant:

\proclaim{Hypothesis C} For all $\ell$, $T_\ell$ is ideally commutative.
\endproclaim

\proclaim{Hypothesis C$'$} For all $\ell$, $T_\ell$ is ideally $\Bbbk$-affine
commutative.
\endproclaim

\subhead 2.11 \endsubhead Assume hypothesis C. 

(i) It immediately follows that $T$ is ideally commutative. 

(ii) If hypothesis A also holds then $Z$ is noetherian, by (i), since $T$ is
right noetherian. If hypothesis C$'$ holds then $Z$ is a finitely generated
$\Bbbk$-algebra.

(iii) As in (2.10), extension and contraction of ideals provide homeomorphisms
between $\Spec T$ and $\Spec Z$, and between $\Max T$ and $\Max Z$.

\subhead 2.12 \endsubhead In the next subsections, before presenting our main
result, we collect information on primitive and related ideals. Our approach
generalizes, and is directly derived from, \cite{\GooLettwo, \S 2}. We first
record the following preparatory lemma:

\proclaim{Lemma} Assume hypotheses {\rm A} and {\rm B}. {\rm (i)} Let
$\p$ be a prime ideal of $T$, and set $P = \eta^*(\p)$. Then $T/\p$ is a 
classical right quotient ring of $R/P$. {\rm (ii)} Every prime factor of $R$ is
right Goldie. \endproclaim

\demo{Proof} (i) Set $Q = \nu^*(\p)$. By (2.7i), $\mu^*(Q) = P$, and by
(2.6ii), $P = Q.T$. Next, it follows from (2.2ii) that $R/P$ is ismorphic to
$S/Q$. As noted in (2.6ii), $T/Q.T = T/\p$ is a classical right quotient ring
of $S/Q$. In particular, part (i) follows.

(ii) First recall that a prime ring is right Goldie if and only if it has a
simple artinian classical right quotient ring; see, for example,
\cite{\GooWar, Chapter 6}, and in particular \cite{\GooWar, 6.18}, for
details. By hypothesis A, $T/Q.T$ is right noetherian and so posesses a simple
artinian classical right quotient ring. However, by (2.6ii), $T/Q.T$ is a
classical right quotient ring for $S/Q$. Therefore, $S/Q$ posesses a simple
artinian classical right quotient ring, and (ii) follows.  \qed\enddemo

\subhead 2.13 \endsubhead We now briefly review the Nullstellensatz, for
noncommutative algebras, and the Dixmier-Moeglin equivalance; the reader is
referred (e.g) to \cite{\McCRob, Chapter 9; \Ren; \Rowone, \S 8.4; \BroGooone,
II.8} for background information.

(i) Using the notation of \cite{\McCRob, 9.1.4}, a $\Bbbk$-algebra $\Lambda$
{\sl satisfies the Nullstellensatz (over $\Bbbk$)\/} if $\Lambda$ satisfies
the {\sl endomorphism property\/} (i.e., $\End M_\Lambda$ is algebraic over
$\Bbbk$ for every simple right $\Lambda$-module $M$), and if $\Lambda$
satisfies the {\sl radical property\/} (i.e., the Jacobson radical of every
factor ring of $\Lambda$ is nil). When $\Bbbk$ is uncountable, every countably
generated $\Bbbk$-algebra satisfies the Nullsetellensatz (cf\., e.g.,
\cite{\McCRob, 9.1.8}). When every prime factor of $\Lambda$ is right or left
Goldie, $\Lambda$ satisfies the radical property if and only if $\Lambda$ is a
Jacobson ring (i.e., every prime ideal is an intersection of primitive
ideals).

(ii) Recall that a prime ideal $P$ of a $\Bbbk$-algebra $\Lambda$ is {\sl
rational (over $\Bbbk$)\/} provided $\Lambda/P$ is right (or left) Goldie, and
provided the center of the Goldie right (or left) quotient ring of $R/P$ is
algebraic over $\Bbbk$.

(iii) Also recall that a prime ideal of a ring $\Lambda$ is {\sl locally
closed\/} if it is a locally closed point of $\Spec \Lambda$. In other words,
$P \in \Spec \Lambda$ is locally closed exactly when
$$Q \; \subsetneq \; \bigcap \{ \, Q' \in \Spec \Lambda \; \mid \; Q'
\supsetneq Q \, \} .$$

(iv) Following \cite{\Ren}, a $\Bbbk$-algebra $\Lambda$ satisfies the {\sl
Dixmier-Moeglin equivalence (over $\Bbbk$)\/} provided a prime ideal $P$ of
$\Lambda$ is primitive if and only if $P$ is rational, if and only if $P$ is
locally closed.

(v) If $\Lambda$ is a Jacobson ring, it follows immediately that every locally
closed prime ideal of $\Lambda$ is primitive. It follows from (i) and (2.12),
when $R$ satisfies the Nullstellensatz, that every locally closed prime ideal
of $R$ is primitive.

(vi) Suppose that $\Lambda$ is a $\Bbbk$-algebra satisfying the
Nullstellensatz, that $P$ is a primitive ideal of $\Lambda$, and that
$\Lambda/P$ is right or left Goldie. A proof that $P$ must then be rational
can be found (e.g.) in \cite{\Dix, 4.1.6}.

\proclaim{2.14 Lemma} Assume hypotheses {\rm A} and {\rm B}, and let $\p$ be a
prime ideal of $T$. Then $\p$ is locally closed if and only if $\eta^*(\p)$ is
a locally closed prime ideal of $R$. \endproclaim

\demo{Proof} It is easy to see that a prime ideal $Q$ of $S$ is locally closed
if and only if $\mu^*(Q)$ is a locally closed prime ideal of $R$. Next, let
$\nu_\ell$ be the canonical map $S_\ell \rightarrow T_\ell$.  Let $\p$ be a
prime ideal of $T_\ell$, and set $Q = \nu^{-1}_\ell(\p)$. Note that $\p =
Q.T_\ell$, as noted in (2.6i). 

To prove the lemma it only remains to check that $\p$ is locally closed in
$\Spec T_\ell$ if and only if $Q$ is locally closed in $\Spec S_\ell$. To
start, when $Q$ is locally closed it follows immediately from (2.6i) that $\p$
must also be locally closed.

So now suppose that $\p$ is locally closed. Let $\Sigma$ be the set of prime
ideals of $S_\ell$ properly containing $Q$, let $\Delta$ be the set of prime
ideals of $T_\ell$ properly containing $\p$, and let $\Gamma$ be the set of
prime ideals $\nu^{-1}(\p')$ for $\p' \in \Delta$. Since the intersection of
the prime ideals in $\Delta$ must strictly contain $\p$, it follows from
(2.6i) that the intersection $K_\ell$ of the prime ideals in $\Gamma$ must
strictly contain $Q$. 

Next, let $L_\ell$ denote the natural image of $J_\ell$ in $S_\ell$. By
(2.6i), $L_\ell$ is not contained in $Q$. It also follows from (2.6i) that if
$Q' \in \Sigma \setminus \Gamma$ then $Q'$ contains $L_\ell$.  We now see that
$$M_\ell \; \colon= \; \bigcap_{Q' \in \Sigma}Q' \; = \; \left(\bigcap_{Q'
\supseteq K_\ell}Q'\right) \bigcap \left(\bigcap_{Q' \supseteq
L_\ell}Q'\right).$$
However, since $Q$ is prime, $K_\ell \cap L_\ell \nsubseteq Q$, and so $M_\ell
\nsubseteq Q$.  Hence $M_\ell$ properly contains $Q$, and $Q$ is locally
closed. The lemma follows.  \qed\enddemo

\subhead 2.15 \endsubhead We now introduce the last hypopthesis:

\proclaim{Hypothesis D} $R$ is a $\Bbbk$-algebra satisfying the
Nullstellensatz. \endproclaim

\proclaim{2.16 Theorem} Assume that $R$ satisfies hypotheses {\rm A}, {\rm B},
{\rm C}, and {\rm D}. Then $R$ satisfies the Dixmier-Moeglin
equivalence. \endproclaim

\demo{Proof} Let $P$ be a rational ideal of $R$ -- recall from (2.12) that
every prime factor of $R$ is right Goldie. By (2.13v, vi), to prove the
theorem it suffices to prove that $P$ is locally closed. Now, by (2.7), we can
choose a prime ideal $\p$ of $T$ such that $\eta^*(\p) = P$. By (2.12i),
$T/\p$ is a classical right quotient ring for $R/P$; hence the rationality of
$P$ ensures the rationality of $\p$. However, since $\p$ is rational, $\p \cap
Z$ is a maximal ideal of $Z$. Therefore, by (2.10) and (2.11i), $\p$
is a maximal ideal of $T$, and in particular, $\p$ is locally closed. It now
follows from (2.14) that $P$ is locally closed.
\qed\enddemo

\subhead 2.17 \endsubhead Choose $\q \in \Spec Z$, and set
$$\theta(\q) \; = \; \ann(T/\q T)_R \; = \; \ann \big( (Z/\q)\otimes
_ZT)_R.$$
By (2.11iii), $\q T$ is a prime ideal of $T$. By (2.7),
$\ann(T/\q T)_R = \eta^*(\q T)$ is a prime ideal of $R$. Consequently,
$\theta(\q)$ is a prime ideal of $R$.

We now present the main result of this section.

\proclaim{2.18 Theorem} Assume hypotheses {\rm A}, {\rm B}, and {\rm C}. 

{\rm (i)} The function $\theta \colon \Spec Z \rightarrow \Spec R$, described
in {\rm (2.17)}, is continuous, finite-to-one, and surjective. If
$X_1,\ldots,X_t$ are disjoint then $\theta$ is bijective.

{\rm (ii)} For all $\q \in \Spec Z$, $\theta(\q)$ is locally closed if and
only if $\q$ is locally closed. 

{\rm (iii)} Assume that $R$ satisfies hypothesis {\rm D}. Then for all $\q \in
\Spec Z$, $\theta(\q)$ is primitive if and only if $\q$ is
maximal. Consequently, $\theta \colon \Max Z \rightarrow \Prim R$ is
continuous, finite-to-one, and surjective.
\endproclaim

\demo{Proof} Part (i) follows from (2.7) and (2.11iii).  Part (ii) follows
from (2.11iii) and (2.14).

To prove part (iii), let $\q$ be a prime ideal of $Z$ and let $P =
\theta(\q)$.  Set $\p = \q T$; by (2.11iii), $\p$ is a prime ideal of $T$, and
following (2.17), $P = \theta(\q) = \eta^*(\p)$. 

Now suppose that $P$ is primitive. By (2.16), $P$ is a rational ideal of
$R$. Therefore, by (2.12i), $\p$ is a rational ideal of $T$, and so $\p \cap
Z$ must be a maximal ideal of $Z$.

Conversely, suppose that $\q$ is maximal.  Then, by (2.11iii), $\p$ is a
maximal ideal of $T$. Therefore, $\theta(\q) = \eta^*(\p)$ is a locally closed
prime ideal of $R$, by (2.14). Thus $\theta(\q)$ is a primitive ideal of $R$,
by (2.16).  Part (iii) follows. \qed\enddemo

\subhead 2.19 \endsubhead Let $\Lambda$ and $\Lambda'$ be rings, and let $U$
be an $\Lambda$-$\Lambda'$-bimodule. In the approach to noncommutative
algebraic geometry developed in \cite{\Ros, \VdB}, the categories $\Mod
\Lambda$ and $\Mod \Lambda'$ are viewed as categories of quasi-coherent
sheaves on not-explicitly-defined noncommutative affine schemes, and the
functor
$$\Mod \Lambda \; @> \quad {\, - \,} \otimes_\Lambda U \quad >> \; \Mod
\Lambda'$$
is an {\sl affine map\/}.

Now assume hypotheses A, B, and C. We obtain the affine map
$$\Mod Z \; @> \quad {\, - \,} \otimes_ZT \quad >> \; \Mod R.$$
Here, of course, $\Mod Z$ is naturally equivalent to the category of
quasi-coherent sheaves on the affine scheme $\Spec Z$, which is noetherian by
(2.11ii). We see that the function $\theta$ of (2.17) is derived from this
affine map.

\subhead 2.20 \endsubhead We conclude this section with useful sufficient
conditions for hypotheses A, B$'$, C, and C$'$. To start, say that a prime
ring $\Lambda$ is {\sl generically ideally commutative\/} provided the right
Ore localization of $\Lambda$ at the non-negative powers of some nonzero
normal (and so regular) element is -- in the notation of (2.10) -- ideally
commutative.  Say that a prime ring $\Lambda$ is {\sl generically ideally
$\Bbbk$-affine commutative\/} provided the right Ore localization of $\Lambda$
at the non-negative powers of some nonzero normal element is $\Bbbk$-affine
ideally commutative.

We will find the following criteria useful:

\proclaim{Proposition} Assume that the semiprime ideals of $R$ satisfy the
ascending chain condition.

{\rm (i)} If every prime factor of $R$ is generically ideally commutative,
then hypotheses {\rm A}, {\rm B$'$}, and {\rm C} are satisfied, for suitable
choices of $I_1,\ldots,I_t$, $J_1,\ldots,J_t$, and $D_1,\ldots,D_t$.

{\rm (ii)} If every prime factor of $R$ is generically ideally $\Bbbk$-affine
commutative, then hypotheses {\rm A}, {\rm B$'$}, and {\rm C$'$} are
satisfied, for suitable choices of $I_1,\ldots,I_t$, $J_1,\ldots,J_t$, and
$D_1,\ldots,D_t$.
\endproclaim

\demo{Proof} (i) Assume every prime factor of $R$ is generically ideally
commutative. Since $\Spec R = V(N)$, where $N$ is the prime radical of $R$, we
may assume without loss of generality that $R$ is semiprime. Moreover, by
noetherian induction, we may further assume that (i) holds for all proper
semiprime factors of $R$.

Suppose that $R$ is not prime. Then there exist nonzero semiprime ideals $I$
and $J$ of $R$ such that $I \cap J = 0$. Since (i) holds for $R/I$ and $R/J$,
and since $\Spec R = V(I) \cup V(J)$, we see that (i) holds for $R$.

So now suppose that $R$ is prime. Since $R$ is generically ideally
commutative, we can choose a regular normal element $b_1$ of $R$ such that the
right Ore localization of $R$ at the non-negative powers of $b_1$ is ideally
commutative. It follows that the conditions in hypotheses A, B$'$, and C hold
for $I_1 = 0$, $J_1 = \langle b_1 \rangle$, $X_1 = V(I_1) \setminus V(J_1)$,
$S_1 = R$, $D_1 = \{1, b_1,b_1^2,\ldots \}$, and $T_1 = S_1D_1^{-1}$.

Now take $I_2 = \sqrt{J_1}$. We know from the induction hypothesis that (i)
holds for $R/I_2$. Since $\Spec R = X_1 \cup V(I_2)$, we see that (i) holds
for $R$.

(ii) This follows similarly to (i). \qed\enddemo

\head 3. PI algebras \endhead 

In this section, assume that $R$ is a ring satisfying a polynomial identity;
see, for example, \cite{\McCRob, Chapter 13} or \cite{\Rowtwo} for necessary
background information.

\subhead 3.1 \endsubhead (i) Let $\Lambda$ be a (prime) Azumaya algebra. In
the notation of (2.10), $\Lambda$ is ideally commutative.  Also, if $\Lambda$
is a finitely generated $\Bbbk$-algebra, then $Z(\Lambda)$ is finitely
generated as a $\Bbbk$-algebra.

(ii) Also in the notation of (2.10), it follows from (i) and the Artin-Procesi
theorem that every prime factor of $R$ is generically ideally
commutative. Moreover, if $R$ is a finitely generated $\Bbbk$-algebra then
every prime factor of $R$ is generically ideally $\Bbbk$-affine commutative.

\subhead 3.2 Finitely generated PI $\Bbbk$-algebras \endsubhead Suppose that
$R$ is a finitely generated $\Bbbk$-algebra. Then the semiprime ideals of $R$
satisfy the ascending chain condition (see, e.g., \cite{\Rowtwo, 4.5.7}). In
particular, by (3.1ii) and (2.20ii), $R$ satisfies hypotheses A, B$'$, and
C$'$. Also, $R$ satisfies hypothesis D (see, e.g., \cite{\McCRob,
13.10.4}). Consequently, for a suitably chosen, finitely generated (by
(2.11ii)), commutative $\Bbbk$-algebra $Z$, (2.18) provides a continuous,
finite-to-one surjection from $\Spec Z$ onto $\Spec R$, restricting to a
continuous surjection from $\Max Z$ onto $\Prim R = \Max R$.

The fact that $\Prim R$ is a disjoint union of finitely many locally closed
subsets, each homeomorphic to an open subset of a $\Bbbk$-affine algebraic
variety, can already be found in \cite{\Art; \Pro}.

\subhead 3.3 Noetherian PI algebras \endsubhead Suppose that $R$ is right
noetherian (but not necessarily finitely generated as an algebra). Then by
(3.1ii) and (2.20i), $R$ satisfies hypotheses A, B$'$, and C. Therefore, (2.18)
provides a continuous, finite-to-one, surjection $\theta\colon \Spec Z
\rightarrow \Spec R$, for a suitably chosen commutative ring $Z$, noetherian
by (2.11ii). If $R$ and $Z$ are Jacobson rings, then, by (2.18ii), $\theta$
restricts to a surjection from $\Max Z$ onto $\Prim R = \Max R$.

\head 4. Solvable Lie algebras \endhead

Assume in this section that $\Bbbk$ has characteristic zero and that $R$ is
the enveloping algebra of a finite dimensional completely solvable $\Bbbk$-Lie
algebra $\frak g$. 

\subhead 4.1 \endsubhead By \cite{\McCRob, 14.9.19}, every prime factor of $R$
is generically ideally $\Bbbk$-affine commutative. Also, $R$ is noetherian by
the Poincar\'e-Birkhoff-Witt theorem. Therefore, by (2.20ii), $R$ satisfies
hypotheses A, B$'$, and C$'$. Furthermore, $R$ satisfies hypothesis D; see,
for example, \cite{\Dix, 2.6.4, 3.1.15}. Hence, (2.18) provides a continuous,
finite-to-one, surjection from $\Spec Z$ onto $\Spec R$, restricting to a
continuous surjection from $\Max Z$ onto $\Prim R$, for a suitable finitely
generated (by (2.11ii)) commutative $\Bbbk$-algebra $Z$.

\subhead 4.2 \endsubhead When $\Bbbk$ is algebraically closed, the Dixmier map
(see, e.g., \cite{\Dix}) provides a Zariski-continuous map from ${\frak g}^*$
onto $\Prim R$; bicontinuity of the factorized map is given in \cite{\Mat}.

\subhead 4.3 \endsubhead An example: Let $\frak g$ be the nonabelian
2-dimensional solvable $\Bbbk$-Lie algebra, and let $R$ be the enveloping
algebra of $\frak g$. Then $R$ is a domain and contains a unique minimal
nonzero prime ideal $P$. Set $I_1 = 0$, $J_1 = P$, $I_2 = P$, and $J_2 =
R$. Following the constructions in (2.2), and in the rest of \S 2, we see that
$Y_1 \cong X_1 = V(0) \setminus V(P) = \{0\}$, and $Y_2 \cong X_2 = V(P)$. Since $0$
is a generic point of $\Spec R$, and since $Y_1$ is closed, we see that
$\theta$ in this case is not a topological quotient map.

\head 5. Quantum Semisimple Groups \endhead 

\subhead 5.1 \endsubhead Let $G$ be a connected, semisimple Lie group,
and let $R$ denote one of the following algebras: $\RqG$ as defined in
\cite{\Josone}, $\CqpG$ as defined in \cite{\HodLevTor}, or $\OqG$ as defined
in \cite{\BroGooone, I.7}. For each of the possible choices of $R$ we let
$\Bbbk$ denote the appropriate ground field (over which $R$ is finitely
generated), and we assume that $\Bbbk$ is algebraically closed. These choices
for $R$ are ``non-root-of-unity'' or ``generic'' quantizations.

Throughout, we equip $R$ with the rational $H$-action by automorphisms
described, for example, in \cite{\BroGooone, II.1.18}, where $H$ is a torus
over $\Bbbk$ corresponding to a maximal torus of $G$. Note, when a
$\Bbbk$-algebra is equipped with an $H$-action by automorphisms, that there is
a naturally induced $H$-action on the prime and primitive spectra.

\subhead 5.2 \endsubhead The following information is well known, and can be
extracted, for example, from \cite{\BroGooone, II.4; \BroGootwo; \HodLevTor,
\S 4; \Josone, Chapter 10}. The original results describing the prime and
primitive ideals in this way can mostly be found in \cite{\HodLevone;
\HodLevtwo; \HodLevTor; \Jostwo}. 

(i) To start, $R$ is a noetherian domain. Also, $R$ satisfies the
Nullstellensatz; see (e.g.) \cite{\BroGooone, II.7.20}.

(ii) Let $W$ be the Weyl group associated to $G$. For each $w \in W{\times}W$
there exists a prime ideal $I_w$ of $R$ and a finite set $E_w$ of elements of
$R$ normal modulo $I_w$ such that
$$\Spec R \; = \coprod_{w \; \in \; W{\times}W} \; \Spec _w R, \qquad
\text{and} \qquad \Prim R \; = \coprod_{w \; \in \; W{\times}W} \; \Prim _w
R,$$
where
$$\Spec_w R \; = \; \{ \,   P \in \Spec R \; : \; \text{$P \supseteq I_w$ and
$P\cap E_w = \emptyset$}  \, \},$$
and
$$\Prim_w \; = \; \big(\Spec_w R\big) \; \cap \; \big(\Prim R\big).$$
Moreover, $\Spec_w R$ and $\Prim_w R$ are nonempty. (That $E_w$ can be chosen
to be finite is explictly justified in \cite{\BroGootwo, 5.6}.)  

Let $b_w$ denote the product of the elements in $E_w$, and set $J_w = \langle
b_w \rangle$. It is easy to see that $b_w$ is normal modulo $I_w$ and that
$$\Spec_w R \; = \; V(I_w) \setminus V(J_w).$$

(iii) Set $S_w = R/I_w$, and let $D_w$ denote the image in $S_w$ of the
non-negative powers of $b_w$. Let $T_w$ denote the localization of $S_w$ at
$D_w$, and let $Z_w$ be the center of $T_w$. Each $T_w$ is ideally
$\Bbbk$-affine commutative, following (e.g.) \cite{\HodLevTor, 4.15}.

(iv) Each $Z_w$ is isomorphic as a $\Bbbk$-algebra to a Laurent polynomial ring
$\Bbbk[y_1^{\pm}, \ldots ,y_{m_w}^{\pm}]$, for some non-negative integer $m_w$
(see, e.g., \cite{\BroGooone, II.4.14}). The value of $m_w$ can be explicitly
determined in one-parameter quantizations and for some multiparameter cases;
see \cite{\HodLevTor, 4.17}.

(v) Under the given action of $H$ on $R$ by automorphisms, each $I_w$ is
$H$-stable and the elements of $E_w$ can be chosen to be
$H$-eigenvectors. There is then an induced $H$-action on $T_w$, and $Z_w$ is
generated by $H$-eigenvectors. The inclusion of $Z_w$ in $T_w$ is
$H$-equivariant, each $\Spec_w R$ is stable under the induced $H$-action, and
each $\Prim_w R$ is a single $H$-orbit.

\subhead 5.3 \endsubhead For the given $I_w$, $J_w$, $b_w$, and $T_w$ of
(5.2), hypotheses A follows from (5.2i), hypothesis B$'$ follows from
(5.2ii), hypothesis C$'$ follows from (5.2iii), and hypothesis D follows
from (5.2i). Moreover, setting $X_w = V(I_w)\setminus V(J_w) = \Spec_w R$,
we see that the $X_w$ are disjoint, as noted in (5.2ii).

Using (2.18) and (5.2v), we obtain an $H$-equivariant, continuous, bijection
$\theta\colon \Spec Z \rightarrow \Spec R$, restricting to an $H$-equivariant
continuous bijection from $\Max Z$ onto $\Prim R$, for $Z$ equal to the ring
direct product of the $Z_w$.

\head 6. Other Quantized Coordinate Algebras \endhead 

In this section, $\Bbbk$ denotes an infinite field, $H$ denotes a
$\Bbbk$-algebraic $n$-torus, and $R$ denotes a finitely generated noetherian
$\Bbbk$-algebra equipped with a rational $H$-action by automorphisms. Assume
further that there are only finitely many fixed points in the naturally
induced $H$-action on $\Spec R$.  In \cite{\BroGooone; \GooLettwo; \GooSta}
this and similar settings are used as unified frameworks to study both the
algebras considered in \S 4 and other quantum function algebras.

\subhead 6.1 \endsubhead That the algebras $\RqG$, $\CqpG$, and $\OqG$
considered in \S 4 satisfy the current assumptions follows from
\cite{\HodLevTor; \Jostwo}; see (5.2). Other examples satisfying the current
assumptions include quantized coordinate rings of affine $n$-space, quantized
Weyl algebras (at non-roots of unity), and quantized coordinate rings (also at
non-roots-of-unity) of $\nxn$ matrices, $GL_n$, symplectic $n$-space, and
euclidean $n$-space; see \cite{\BroGooone, II; \GooLettwo} for details. All of
these examples satisfy hypothesis D; see (e.g.) \cite{\BroGooone, II.7.17,
II.7.20}.

\subhead 6.2 \endsubhead We now verify hypotheses A, B, and C$'$, for
suitable choices of $I_\ell$, $J_\ell$, and $D_\ell$. To start, let
$I_1,\ldots,I_t$ be the $H$-stable prime ideals of $R$, and for each $1 \leq
\ell \leq t$, let $J_\ell$ be the intersection of the $H$-stable prime ideals
of $R$ not contained in $I_\ell$. 

(i) Following \cite{\BroGooone, II.2; \GooLettwo, \S 2; \GooSta},
$$\Spec R \; = \; \coprod _{\ell = 1}^t X_\ell,$$
where
$$X_\ell \; = \; V(I_\ell) \setminus V(J_\ell),$$
for all $1 \leq \ell \leq t$. 

(ii) For each $\ell$, set $S_\ell = R/I_\ell$. It follows from \cite{\GooSta}
(or \cite{\BroGooone, II.2.13}) that there exist right denominator sets
$D_\ell$ of $S_\ell$, for all $\ell$, such that hypothesis B holds true.  Now
set $T_\ell = S_\ell D_\ell^{-1}$. Since $R$ is noetherian, each $T_\ell$ is
noetherian, and so hypothesis A is satisfied. Let $Z_\ell$ denote the center
of $T_\ell$, for each $\ell$. Hypothesis C$'$ holds in this setting, by
\cite{\GooSta} (or \cite{\BroGooone, II.2.13}). 
 
\subhead 6.3 \endsubhead Retain the notation of (6.2). 

(i) Since each $I_\ell$ is $H$-stable, it follows that each $S_\ell$ inherits
the $H$-action on $R$. By \cite{\GooSta} (or \cite{\BroGooone, II.2.13}), we
can assume that $D_\ell$ consists entirely of $H$-eigenvectors. Hence each
$T_\ell$ is naturally equipped with a rational $H$-action by automorphisms.
Furthermore, it follows from \cite{\GooSta} (or see \cite{\BroGooone,
II.2.13}) that each $Z_\ell$ is generated by $H$-eigenvectors.

(ii) Set $T = T_1 \times \cdots \times T_t$, and identify the center $Z$ of $T$
with $Z_1 \times \cdots \times Z_t$. Then $T$ and $Z$ inherit rational
$H$-actions by automorphisms, from the $H$-actions on each $T_\ell$ and
$Z_\ell$. The embedding of $Z$ in $T$ is $H$-equivariant, as is the
$\Bbbk$-algebra homomorphism $\eta\colon R \rightarrow T$ of (2.4). Following
(2.11ii), $Z$ is finitely generated as a $\Bbbk$-algebra.

\subhead 6.4 \endsubhead By (6.2), hypotheses A, B, and C$'$ hold for the
given choices of $I_\ell$, $J_\ell$, and $D_\ell$; moreover, the sets $X_\ell$
are (pairwise) disjoint. In view of (6.3) and (2.18), we obtain a continuous,
$H$-equivariant bijection $\theta\colon \Spec Z \rightarrow \Spec R$.

All of the algebras mentioned in (6.1) satisfy the Nullstellensatz over
$\Bbbk$, by \cite{\BroGooone, II.7.18, II.7.20}; in this case, by (2.18), the
preceding function $\theta$ restricts to a continuous, $H$-equivariant
bijection from $\Max Z$ onto $\Prim R$.

\subhead 6.5 \endsubhead Suppose that $R$ is the quantized coordinate ring of
$\Bbbk$-affine $n$-space, and assume further that $-1$ is not contained in the
multiplicative group generated by the quantizing parameters. It is proved in
\cite{\GooLetone} that $\Spec R$ is a topological quotient of $\Spec
\Bbbk[x_1,\ldots,x_n]$ and that $\Prim R$ is a topological quotient of $\Max
\Bbbk[x_1,\ldots,x_n]$. Goodearl has conjectured in \cite{\Goo} that the prime
spectra of more general quantum function algebras should be topological
quotients of commutative prime spectra.

\enddocument